\documentclass[12pt,reqno]{amsart}

\usepackage{enumitem}

\setlist[enumerate]{itemsep=1ex}
\setlist[itemize]{itemsep=1ex}

\usepackage{amsmath,amssymb,amsthm,tikz-cd,psfrag,fourier}
\usepackage{microtype}
\usepackage[a4paper, total={6.5in, 8.5in}]{geometry}

\usepackage{graphicx}

\usepackage[final,
pdftitle={TITLE},
pdfauthor={NAME},hypertexnames=false]{hyperref}

\usepackage{epstopdf}

\theoremstyle{plain}

\newtheorem*{rem*}{Remark}
\theoremstyle{definition}
\newtheoremstyle{case}{}{}{}{}{}{:}{ }{}

\numberwithin{subcase}{case}

\theoremstyle{plain}
\newtheorem{thm}{Theorem}

\newtheorem{que}[thm]{Question}

\theoremstyle{definition}

\newtheorem{remark}[thm]{Remark}

\def\scl[#1][#2]{{\scalebox{#1}{#2}}}

\def\stab[#1]{\mathrm{Stab}(#1)}

\def\cob[#1]{\textcolor{blue}{#1}}
\def\cog[#1]{\textcolor{green}{#1}}
\def\cor[#1]{\textcolor{red}{#1}}

\newcommand{\param}%
	{{\mathchoice{\mkern1mu\mbox{\raise2.2pt\hbox{$\centerdot$}}\mkern1mu}%
	{\mkern1mu\mbox{\raise2.2pt\hbox{$\centerdot$}}\mkern1mu}%
	{\mkern1.5mu\centerdot\mkern1.5mu}{\mkern1.5mu\centerdot\mkern1.5mu}}}

\begin{document}


\title[SA SFTs on hyperbolic groups]{Strongly aperiodic SFTs on hyperbolic groups: \\ where to find them and why we love them}

\author[Yo'av Rieck]{Yo'av Rieck}
\address{Dept.\ of Mathematics \\
  University of Arkansas\\
  Fayetteville, AR 72701}
\email{\href{mailto:yoav@uark.edu}{yoav@uark.edu}}

\begin{abstract}
In~\cite{cohen_goodman-strauss_rieck_2021} D. B. Cohen, C. Goodman--Strauss, and the author  proved that a hyperbolic group admits an ``SA SFT'' if and only if it has at most one end. This paper has two distinct parts: the first is a conversation explaining what an SA SFT is and how they may be of use. In the second part I attempt to explain  both old and new ideas that go into the proof. References to specific claims in~\cite{cohen_goodman-strauss_rieck_2021} are given, with the hope that any interested reader may be able to find the details there more accessible after reading this exposition.
\end{abstract}

\clearpage\maketitle
\thispagestyle{empty}

%
%


This paper is about a structure abbreviated as SA SFT, which stands for strongly aperiodic subshift of finite type, that exists on some groups and not on others (all groups considered are infinite and finitely generated).  It has two parts: in the Part~\ref{part:why} we discuss the reasons for constructing SA SFTs and explain what they are. In Part~\ref{part:how} we explain the basic ideas behind the construction which led to the theorem of D. B. Cohen, C. Goodman--Strauss, and the author~\cite{cohen_goodman-strauss_rieck_2021}:

\begin{thm}
\label{MainThm}
A hyperbolic group admits an SA SFT if and only if it has at most one end.     
\end{thm}

\noindent
This paper is written for topologists and geometric group theorists, and we are not assuming any familiarity with SA SFTs.  We do assume familiarity with hyperbolic groups, that were introduced by Gromov in an extremely influential paper~\cite{gromov}. Gromov outlined the study of SFTs on hyperbolic groups in that paper, and a detailed investigation was carried out by Coornaert and Papadopoulos~\cite{MR1222644,MR1878587}. Our interest is promoting the SFT to an SA SFT.

\begin{remark}
Much of this paper reviews old and well-known results. In order to facilitate the reading of~\cite{cohen_goodman-strauss_rieck_2021}, we included many references to that work. This is \em not \em an indication that all the claims are original. 
\end{remark}

{\bf Acknowledgement.}

I would like to thank the organizers of the special session, Jos\'e Ayala Hoffman, Mario Eudave Mu\~noz, and Jennifer Schultens, for giving me this opportunity. I thank Chaim Goodman--Strauss for help preparing this paper and the beautiful illustrations.  
I am very grateful to the anonymous referee for very helpful suggestions.
Yo'av Rieck was partially supported by the Simons Foundation (award number 637880).

\part{Why do we do what we do?}\label{part:why}

\section{The Geography of a Group}

Fix, for the entirety of Part~1, an infinite, finitely generated group \(G\), together with a finite generating set \(S = S^{-1}\). 

As usual we view the group as its Cayley graph \(\Gamma\) which we define here for completeness.  The vertices of \(\Gamma\) are the elements of \(G\), and the edges out of a vertex \(g\) have the form \((\textit{g},\textit{ga})\) for all \(\textit{a} \in S\). The edges of \(\Gamma\) are labeled and directed --- the edge \((\textit{g,ga})\) is labeled \(a\) and points from \(g\) to \(\textit{ga}\). Since \(S = S^{-1}\) there is an edge pointing in the opposite direction, namely, \((\textit{ga},\textit{gaa}^{-1})\).  Connectivity of \(\Gamma\) follows from the fact that \(S\) generates \(G\). Assigning length \(1\) to each edge induces a metric on \(\Gamma\). Left multiplication defines an action of \(G\) on \(\Gamma\) by a graph isomorphism that preserve the labels, directions, and distances. This action is transitive on the vertices.

Since \(G\) acts on the vertices of \(\Gamma\) transitively, we cannot use \(\Gamma\) to distinguish elements, that is, for any \(g_{1},g_{2} \in G\) and any \(R>0\) the balls \(B(R,g_{1})\) and  \(B(R,g_{2})\) are identical. 

How are we to know where we are?

One way is to name each and every group element. For example, if the group is the integers \(\mathbb{Z}\), each integer has a name and we may use it. This will not do. When labeling group elements we may only use \em finitely many \em labels (imagine a small device attached to each group element that records the label) and the group \(G\) is, by assumption, infinite.

Here's another suggestion: label the identity by \(0\) and all other elements \(1\). This poses a different problem: once the labeling of the group is completed we want to be able to \em guarantee \em  that it is ``correct'', and this must be done locally. If we allow arbitrarily large balls to be labeled \(1\), what is to stop us from labeling the entire group \(1\)?

This pinpoints our goal: to construct a finite set of labels and a collection of local rules that will allow us to distinguish elements of \(G\). Specifically, we would like to (and here we give an intuitive description of an SA SFT):
\begin{enumerate}
\item Choose a finite set \(A\).
\item Construct \em local rules \em for labeling.
\item In \em any \em labeling that obeys the local rules, and for any \(g_{1} \neq g_{2} \in G\), there is \(R>0\), so that the labeling on  \(B(R,g_{1})\) and  \(B(R,g_{1})\) are distinct (\(R\) depends on \(g_{1}\) and \(g_{2}\) and this is unavoidable).
\item There must be some labeling obeying the rules, so condition~(3) is not vacuous. 
\end{enumerate}

\begin{remark}
\label{remark:SFT}
A few remarks are in order:
\begin{enumerate}
\item A finite set carries no information other than its cardinality. (In reality, in order to work with \(A\), each label will carry useful information; a similar comment can be made about the states of a Turing machine.)
\item The local rules are the key here. \em Local rules \em form an atlas, or a finite collection of {\it charts}, where each chart is a labeling of a ball of radius \(R\), for some fixed \(R\).\footnote{  
The reader may have seen local rules defined on finite sets that are not necessarily balls, but there is no loss of generality in assuming that they are all balls and of the same radius.}
We say that the labeling satisfies the local rules if the labeling about each point \(g\) agrees with one of the charts (translated to be centered at \(g\)).
\item From this point on we focus on labeling the vertices, but the fixed labelling of the edges plays a key role: since the edges are labeled and directed, once we center a chart at \(g\) the exact position of each vertex of the ball is determined (we cannot ``rotate'' the ball). 

\item An SFT is only required to satisfy the first two conditions above: a finite set, and local rules for labeling. The other two conditions may not be satisfied: on the one hand, the rules may allow for a labeling in which some (even all!) vertices look the same; on the other hand,  the rules may not allow for \em any \em labeling at all (this is an empty SFT). Adding conditions~(3) and~(4) is the content of ``SA''.

\end{enumerate}
\end{remark}

\section{Relation to Tilings}
\label{tiling}

An SFT on a group \(G\) is sometimes called a ``tiling'' of \(G\). This is very reasonable, as we now explain. Starting from the second point of Remark~\ref{remark:SFT}, we see that an SFT is given by an atlas of charts. Say there are \(n\) charts and number them \(1\) through \(n\).

Now the labelling about each \(g\), with labels from \(A\), must coincide with one of the charts. We note that, by writing the number of the chart at \(g\) (thus obtaining a function \(G \to \left\{1,\dots, n \right\}\)). Conversely, any function \(G \to \left\{1,\dots, n\right\}\) gives a labeling \(G \to A\)  provided that overlapping charts are compatible.  An easy argument shows that if the labeling about neighboring vertices agree, then the labeling is consistent.  

And so here is the tiling: each tile covers exactly one group element \(g\) and carries a label from \(\{1,\dots,n\}\) that corresponds to the chart used for the labeling about \(g\). Let \(T_{g}\) denote the tile that covers \(g\).  
We can visualize \(T_{g}\) as a ``polygon'' that covers only \(g\), where the boundary of \(T_{g}\) consists of finitely many ``edges'', one for each neighbor of \(g\), where ``edges'' can be glued together if the corresponding charts are compatible.  Here ``edge'' means an ``edge'' of the ``polygon'', not an edge of the Cayley graph (in fact every ``edge'' of the glued ``polygons'' crosses exactly one edge of the Cayley graph and {\it vice-versa}). 
By construction, \(T_{g}\) depends only on the chart at \(g\); thus there are exactly \(n\) possible tile shapes. For obvious reasons this is called  a \em nearest neighbor SFT\em.

We conclude that for a group \(G\) the concepts of \em SFT \em and \em tiling\em\ are the same. The notions of ``aperiodicity'' (to be defined in the next section) carry over. We will focus on SFTs.

\newpage

One final remark: SFTs on \(\mathbb{Z}^{2}\) correspond to tilings of \(\mathbb{E}^{2}\) (where by \(\mathbb{E}^{2}\) we mean \(\mathbb{R}^{2}\) with the Euclidean metric), but the converse does not hold. There are tilings of \(\mathbb{E}^{2}\) (for example~\cite{radin}) that are not aligned along a lattice and hence do not induce an SFT on \(\mathbb{Z}^{2}\). Similarly, strongly aperiodic tiles for \(\mathbb{H}^{2}\) were first constructed by Goodman-Strauss in~\cite{MR2142334}, but do not provide an SFT on any group of isometries of \(\mathbb{H}^{2}\); the first SFT on a hyperbolic was constructed by Cohen and Goodman--Strauss~\cite{MR3692905} 12 years later.

\section{Formal Definitions}

After reviewing the concepts the formal definitions should be easy to read. For a finitely generated group \(G\) we define:

\begin{enumerate}
\item The full shift on \(G\) with labels in a finite set \(A\) is: 
\[
A^{G} := \left\{ \omega : G \to A\right\}
\]
In other words, these are all possible labeling of the elements of \(G\) by labels from \(A\).

\item We endow \(A\) with the discrete topology and \(A^{G}\) with the product topology. This turns \(A^{G}\) into a compact space that carries an right \(G\)-action; the action of \(g \in G\) on \(\omega \in A^{G}\) is defined as follows (here, we evaluate \(\omega \cdot g \in A^{G}\) on \(h \in G\)):
\[
(\omega \cdot g) (h) = \omega(gh)
\]
This action is by homeomorphisms.

\item A \em subshift \em is a closed, invariant subset \(\Omega \subset A^{G}\). The action of \(G\) on \(A^{G}\) induces an action of \(G\) by homeomorphisms on the compact space \(\Omega\). Note that \(\Omega = \emptyset\) is a subshift. More interesting examples are given by the closure of the orbit of any \(\omega \in A^{G}\).

\item Given a subshift \(\Omega\) we call \(\omega \in \Omega\) a \em configuration\em. Thus any configuration is a function
\[
\omega: G \to A
\]
Consistent with the terminology above, we call the value \(\omega(g)\) the \em label \em at \(g\).

\item A subshift is called \cor[\em strongly aperiodic (SA)] if it is not empty, and for every \(\omega \in \Omega\) we have that \textcolor{red}{\(\stab[\omega]\) is trivial} (here and below \(\stab[\cdot]\) stands for \em stabilizer\em).

\item A non-empty subshift is called \em weakly aperiodic  \em  if \(\stab[\omega]\) has infinite index in \(G\) for every \(\omega \in \Omega\). Since we restricted ourselves to infinite groups, strong aperiodicity implies weak aperiodicity (but at this point it is not clear if the converse holds). We will focus on strong aperiodicity. 

\item A \em subshift of finite type (SFT) \em is a subset of the full shift \(A^{G}\) that is defined by local rules, or an atlas of charts, as described above. It is not defined as a subshift satisfying an extra condition, but in fact it is: it is not hard to see that an SFT is closed and \(G\)-invariant, and hence a subshift.

\item Combining~(5) and~(7), we define {\it a strongly aperiodic subshift of finite type (SA SFT)} to be an SFT in which the stabilizer is trivial for every configuration.

\end{enumerate}

\section{Origins}
The origins of the theory, in the very early 60's, are closely linked to axiomatization and computational theory (as we shall see shortly). H. Wang~\cite{MR1112395} asked if the {\it Domino Problem} is decidable for tilings of \(\mathbb{E}^{2}\): 

{\bf The Domino Problem.} Is there an algorithm that decides if a given finite set of tiles in \(\mathbb{E}^{2}\) can be used to tile the plane? 

Obviously, one can replace \(\mathbb{E}^{2}\) with other spaces of interest, for example, \(\mathbb{H}^{2}\), \(\mathbb{E}^{n}\), and \(\mathbb{H}^{n}\). (We will not discuss higher dimensions.)

Wang proposed an ``algorithm'': try! If the tiles do not tile \(\mathbb{E}^{2}\), an attempt to tile arbitrarily large balls would eventually fail (and one will certainly detect this). On the other hand, if they do tile, one ``should'' find a tiled domain. Of course, in the latter the tiling would be a lift of a tiling of the torus. But does this actually work? Perhaps there is a set of tiles that tiles the plane but not the torus? This seemed rather unlikely. 

\begin{figure}
\centerline{{\includegraphics[]{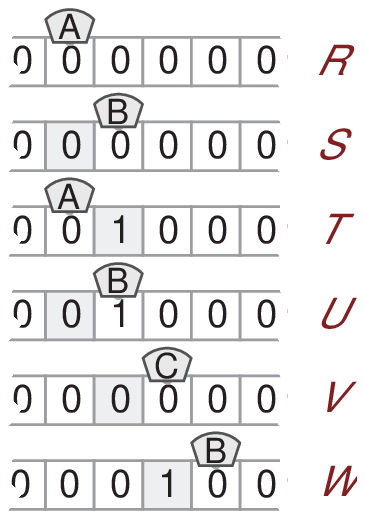}} \hspace{1in}{\includegraphics[]{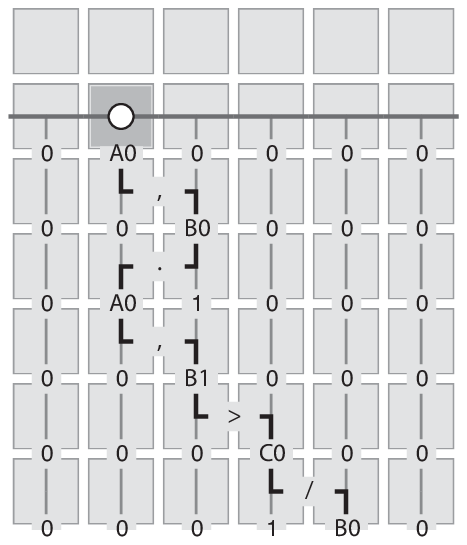}}}
\caption{Wang emulates the run of a Turing machine as a tiling problem.}
\label{wangrun}
\end{figure}

Wang himself found such tiles, but they are \em seeded\em, which means that they contain a special tile (the seed) that must be used. With this, Wang constructed tiles that emulate an arbitrary Turing machine and tile the plane if and only if the given machine never halts; undecidability of the Domino Problem for seeded tiles follows from undecidability of the Halting Problem.  Figure~\ref{wangrun} indicates a Turing machine and seeded tiles that emulate it. Note that without the seed (in the second row, second column) copies of the blank tiles (in the first row) can always be used to tile \(\mathbb{E}^{2}\) (even the torus), and thus without the seed these tiles give nothing of interest.

\psfrag{R}[ll][cc][.8]{$t=0$}
\psfrag{S}[ll][cc][.8]{$t=1$}
\psfrag{T}[ll][cc][.8]{$t=2$}
\psfrag{U}[ll][cc][.8]{$t=3$}
\psfrag{V}[ll][cc][.8]{$t=4$}
\psfrag{W}[ll][cc][.8]{$t=5$}
\psfrag{X}[ll][cc][.8]{$t=6$}  
\psfrag{1}[cc][cc][.8]{${\tt 1}$}    
\psfrag{0}[cc][cc][.8]{${\tt 0}$}
\psfrag{A}[cc][cc][.8]{${\tt A}$}
\psfrag{B}[cc][cc][.8]{${\tt B}$}
\psfrag{C}[cc][cc][.8]{${\tt C}$}
\psfrag{A0}[cc][cc][.8]{${\tt A  0}$}
\psfrag{B0}[cc][cc][.8]{${\tt B 0}$}
\psfrag{C0}[cc][cc][.8]{${\tt C 0}$}
\psfrag{A1}[cc][cc][.8]{${\tt A 1}$}
\psfrag{B1}[cc][cc][.8]{${\tt B 1}$}
\psfrag{C1}[cc][cc][.8]{${\tt C 1}$}
\psfrag{pA0}[cc][cc][.6]{$\phi({\tt A  0})$}
\psfrag{pB0}[cc][cc][.6]{$\phi({\tt B 0})$}
\psfrag{pC0}[cc][cc][.6]{$\phi({\tt C 0})$}
\psfrag{pA1}[cc][cc][.6]{$\phi({\tt A 1})$}
\psfrag{pB1}[cc][cc][.6]{$\phi({\tt B 1})$}
\psfrag{pC1}[cc][cc][.6]{$\phi({\tt C 1})$}
\psfrag{A0}[cc][cc][.8]{${\tt A  0}$}
\psfrag{B0}[cc][cc][.8]{${\tt B 0}$}
\psfrag{C0}[cc][cc][.8]{${\tt C 0}$}
\psfrag{A1}[cc][cc][.8]{${\tt A 1}$}
\psfrag{B1}[cc][cc][.8]{${\tt B 1}$}
\psfrag{C1}[cc][cc][.8]{${\tt C 1}$}
\psfrag{H}[cc][cc][.8]{${\tt H}$}
\psfrag{!}[cc][cc][.6]{$\phi({\tt A  0})$}
\psfrag{@}[cc][cc][.6]{$\phi({\tt B 0})$}
\psfrag{#}[cc][cc][.6]{$\phi({\tt C 0})$}
\psfrag{^}[cc][cc][.6]{$\phi({\tt A 1})$}
\psfrag{&}[cc][cc][.6]{$\phi({\tt B 1})$}
\psfrag{*}[cc][cc][.6]{$\phi({\tt C 1})$}
\psfrag{,}[cc][cc][.6]{$\phi({\tt A  0})$}
\psfrag{.}[cc][cc][.6]{$\phi({\tt B 0})$}
\psfrag{/}[cc][cc][.6]{$\phi({\tt C 0})$}
\psfrag{<}[cc][cc][.6]{$\phi({\tt A 1})$}
\psfrag{>}[cc][cc][.6]{$\phi({\tt B 1})$}
\psfrag{?}[cc][cc][.6]{$\phi({\tt C 1})$}

\begin{figure}
\centerline{{\includegraphics{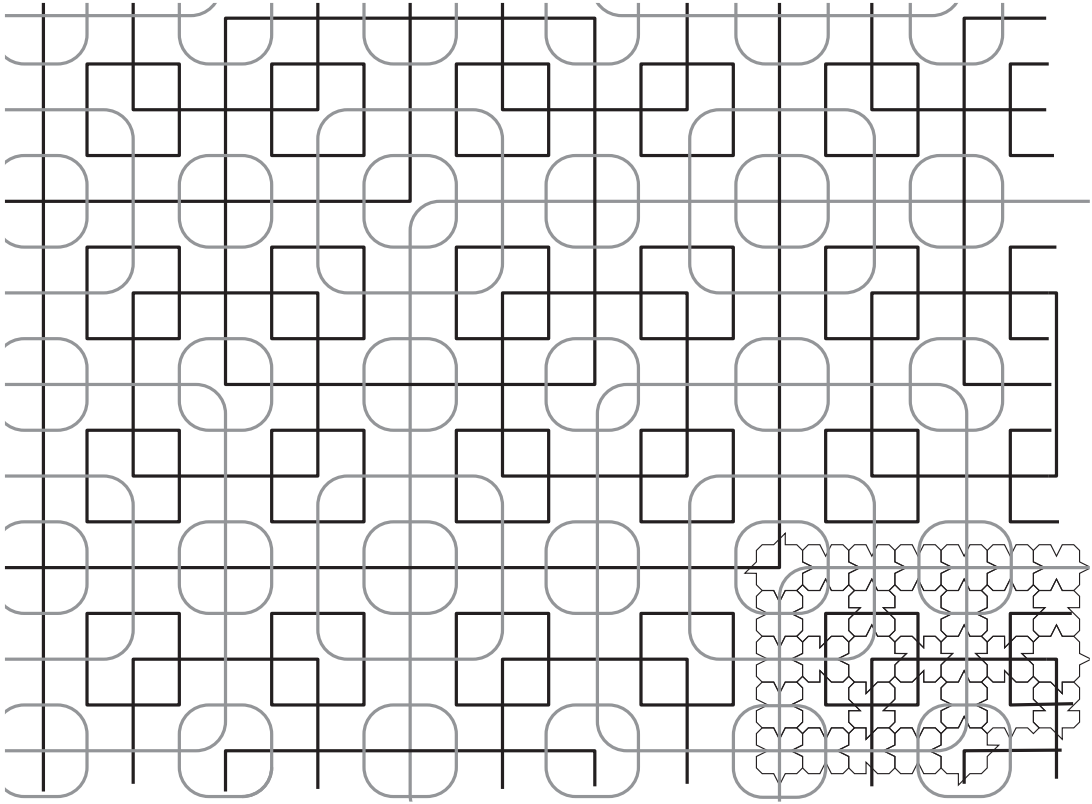}}}
\caption{Unseeded tiles.}
\label{UnseededTiles}
\end{figure}

A few years later, R. Berger showed how to emulate an arbitrary Turing machine with an {\it unseeded} set of tiles (see Figure~\ref{UnseededTiles}, which shows a much-simplified set of tiles due to R. Robinson).  These tiles can be used to show two things:
\begin{enumerate}
\item The Domino Problem in \(\mathbb{E}^{2}\) is undecidable. 
\item There exists a set of strongly aperiodic tiles for \(\mathbb{E}^{2}\). This is indicated in Figure~\ref{UnseededTiles}. As suggested by that figure, the tiles assemble to form disjoint squares of arbitrarily large scale. Hence any element of the stabilizer will need to move these square off themselves or not move them at all, showing that the stabilizer is trivial.
\end{enumerate}


Later, J. Kari~\cite{MR1417578} gave a different and very simple construction that also implies undecidability of the Domino Problem in \(\mathbb{E}^{2}\).

\section{Other Geometries}

Existence of SA SFT's has subtle connections to the geometry of the group. For example, D. Cohen~\cite{DBCohen} showed that for finitely presented groups, existence is a quasi-isomotery invariant.\footnote{For the definition of quasi-isomotery see, for example,~\cite{BridsonHaefliger}.} 
I find this quite surprising because a quasi-isometry is very coarse and can destroy any local information, and SFT's are defined using local rules only. 

Here is a sample of results. Some show existence of an SA SFT, and some show that the Domino Problem is undecidable (implying the existence of a weakly aperiodic set of tiles, but not implying the existence of an SA SFT). 
More results can be found in~\cite{cohen_goodman-strauss_rieck_2021}. The 2018 survey~\cite{DominoSurvey} contains many results about the Domino Problem.  All groups considered in this section are assumed to have decidable word problem and to be finitely presented (although some of the results require only finite generation). 
 
\begin{figure}
\centerline{{\includegraphics[width=4.5in]{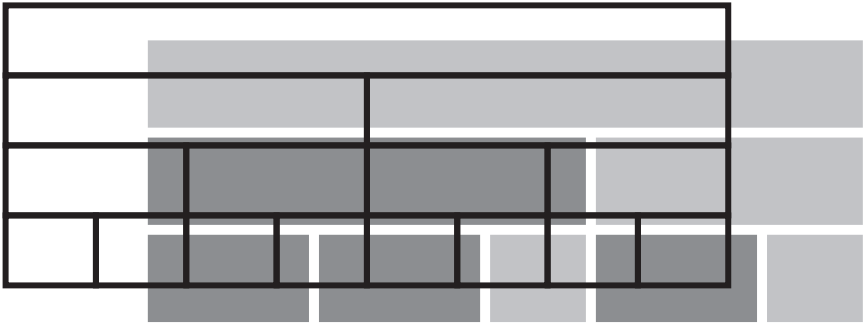}}}
\caption{Two patterns of ``rectangles'' are shown, each rectangle having some predecessor above
and some successors below. In the pattern drawn with dark lines, the number of rectangles
doubles from row to row. In the gray pattern, light rectangles (which are all congruent)
have one light and one dark rectangle as successors, and dark rectangles (which are all
congruent) have one light and two dark successors.}
\label{SA_SFT_H2}
\end{figure}

In 1997 Mozes~\cite{mozes} constructed strongly aperiodic tilings on any simple Lie group \(\Gamma\) of rank greater than one. The tiles themselves are Voronoi cells of a uniform lattice \(G \leq \Gamma\).\footnote{The voronoi cell  corresponding to \(g \in G\) consists the all \(x \in \Gamma\) satisfying \(d(x,g) \leq d(x,g')\) for all \(g' \in G\).}  
The Voronoi cells are then decorated with \em combinatorial \em information, given by the pattern of intersection they form with the Voronoi cells of a second, incompatible, uniform lattice. Mozes showed that incompatibility of the lattices breaks all symmetry and yields a strongly aperiodic tiling of \(\Gamma\); for our purposes, it also gives an SA SFT on \(G\).
(This situation is similar to the tiling shown in Figure~\ref{SA_SFT_H2}, although these tiles come from a different construction in \(\mathbb{H}^{2}\).)

We remark that uniform lattices in rank-one simple Lie groups are hyperbolic and therefore admit an SA SFT by Theorem~\ref{MainThm}. Thus a uniform lattice in any simple Lie group admits an SA SFT.

In~\cite{MR2142334} Goodman-Strauss constructed the first set of strongly aperiodic tiles in \(\mathbb{H}^{2}\). A simpler construction is shown in Figure~\ref{SA_SFT_H2} where the grid ``rectangles'' are all isometric pentagons, each with five edges of length \(1\), fours vertices of angle \(\pi/2\), and one vertex of angle \(\pi\), in the middle of the bottom edge. These overlay an incompatible shaded grid that has two types of tiles, hexagons and pentagons. We note however that these tiles do not correspond to a lattice and therefore do not provide an example of a SA SFT.

In a recent preprint~\cite{ProdusctOfNonAmenable}, Barbieri, M. Sablik, and V. Salo showed that the product two non-amenable groups with decidable word problem admits an SA SFT. 

In~2019 S. Barbieri~\cite{MR3964145} showed that the product of three infinite groups, each with decidable word problem, admits an SA SFT.

In~\cite{MR3594268} Aubrun and Kari constructed SFTs on   \(\mathrm{BS}(1,n)\) groups, that were shown to be SA  by J. Esnay and E. Moutot in~\cite{esnay2021weakly}. Furthermore, in a recent paper~\cite{AUBRUN2021}, Aubrun and Kari showed that the Domino Problem is undecidable for all \(\mathrm{BS}(m,n)\). 

The reader may have noticed that not even one of the groups discussed above is hyperbolic. This changed in 2017 when Cohen and Goodman--Strauss showed that surface groups admit an SA SFT~\cite{MR3692905}. Since then Aubrun, Barbieri, and Moutot~\cite{DominoForSurfaceGroups} showed that the Domino Problem for surface groups is undecidable.

\section{Two Necessary Conditions}
\label{section:TwoNecessaryConditions}

From this point on we only considering finitely presented groups.
There are two known necessary conditions for existence of an SA SFT on a group \(G\):

	\begin{enumerate}
	\item \(G\) must have a decidable word problem.
	\item \(G\) must have at most one end.
	\end{enumerate}
	
This leads us to ask:

\begin{que}
Does every finitely presented, one-ended group with decidable word problem admit an SA SFT?
\end{que}

The first condition, proved by E. Jeandel~\cite{jeandel}, fits well with the philosophy that an SA SFT gives us a way to ``address'' group elements. After all, being able to construct arbitrarily large balls of the Cayley graph is equivalent to having a decidable word problem. 

Here are a few details regarding the second condition.

It was shown by Cohen~\cite{DBCohen} that no group with more than one end admits an SA SFT. We demonstrate this for \(\mathbb{Z}\), which is a well-known and easy case.  Suppose we are given a non-empty SFT on \(\mathbb{Z}\) with label set \(A\) and local rules that are defined on sets of (say) \(n\) adjacent integers. Note that an interval of length \(n\) can be labeled in finitely many ways. It follows that there are two disjoint intervals of length \(n\) with identical labeling; say that one starts at \(a\) and the other starts at \(b > a\). It is now an easy exercise to show that the labeling of the integers from \(a\) to \(b-1\) can be repeated to produce a periodic configuration.

Let us return for a moment to Wang's tiles. Recall that these tiles had a \em seed\em\ that we were required to use. The following example demonstrates the dramatic effect of seeding: 
\begin{itemize}
\item \(A := \{\mathbf{C},\mathbf{L},\mathbf{R}\}\). The meaning of the labels: \(\mathbf{C}\) is the ``center'' of \(\mathbb{Z}\), \(\mathbf{L}\) and \(\mathbf{R}\) mean ``left/right of \(\mathbf{C}\)''.
\item Local rules (allowable configurations, defined on adjacent pairs):
	\begin{itemize}
	\item \(\mathbf{C}\)-\(\mathbf{R}\)
	\item \(\mathbf{L}\)-\(\mathbf{C}\)
	\item \(\mathbf{R}\)-\(\mathbf{R}\)
	\item \(\mathbf{L}\)-\(\mathbf{L}\)
	\end{itemize}
\item The seed is \(\mathbf{C}\). 
\end{itemize}
A moment reflection shows that any allowable configuration has the form
\[
\cdots\textrm{-}
\mathbf{L}\textrm{-}\cdots\textrm{-}\mathbf{L}
\textrm{-}\mathbf{C}\textrm{-}
\mathbf{R}\textrm{-}\cdots\textrm{-}\mathbf{R}
\textrm{-}\cdots
\]
Of course, the stabilizer of such a configuration is trivial, since any non-trivial translation will move \(\mathbf{C}\) to a point labeled \(\mathbf{L}\) or \(\mathbf{R}\). So the seeded SFT is SA. This is in contrast to \em unseeded \em SFT's, where we saw that \(\mathbf{Z}\) admits no SA SFT. Indeed, if we remove the seeding requirement we have two periodic configurations:
\[
\cdots\textrm{-}
\mathbf{R}\textrm{-}\cdots\textrm{-}\mathbf{R}
\textrm{-}\cdots
\ \ \ \text{and} \ \ \
\cdots\textrm{-}
\mathbf{L} \textrm{-} \cdots \textrm{-} \mathbf{L}
\textrm{-} \cdots
\]

This concludes the first part, in which we attempted to do three things: explain what an SA SFT is, motivate our interest in them, and survey some of the many results. The reader should be aware that in recent years there has been a flurry of activity and this discussion is far from complete.

\part{How do we do what we do?}\label{part:how}
In the second part of this paper we discuss the construction of an SA SFT on a 1-ended hyperbolic group \(G\) which we fix once and for all. We also fix a finite generating set \(S = S^{-1}\).  We start with the construction of an SFT that is not required to be strongly aperiodic, and has been known for a long time (\(\Omega_{S}\) below). We then enhance it to attain the desired SA SFT. 

A finitely presented group \(G\) (with a fixed finite generating set \(S = S^{-1}\)) is called {\it hyperbolic} if there is some \(\delta>0\) so that every triangle in the Cayley graph of \(G\) is \(\delta\)-slim, that is, any point on one edge is within \(\delta\) from the union of the other two edges. We fix \(\delta>0\) for the remainder of this paper.

Hyperbolicity, due to Gromov~\cite{gromov}, turns out to be both very natural (satisfied by many groups that appear in application) and very useful. As a simple example, hyperbolicity implies a very efficient solution to the word problem known as \em Dehn's algorithm\em; surprisingly, this turns out to be equivalent to hyperbolicity.  As a far more sophisticated example we mention that the isomorphism problem for hyperbolic groups is decidable (Z. Sela for torsion free hyperbolic group~\cite{MR1324134}, later extended to all hyperbolic groups by F. Dahmani and V. Guirardel~\cite{MR2795509}). These are, of course, just examples of what one can do with hyperbolic groups; very many other beautiful results are known and at this point the theory of hyperbolic groups seems very well understood.

An excellent reference is~\cite{BridsonHaefliger}. Another very useful reference is~\cite{Word_processing_in_groups}, in particular for the {\it shortlex FSA} which will be use extensively below.

\section{Shortlex shellings}
\label{section:Shortlex}

In this section we explain how to construct a certain SFT, called the \em shortlex SFT \em and denoted \(\Omega_{S}\), on a hyperbolic group \(G\); these SFT's will serve as the backbone for our main construction. The group \(G\) need not be one-ended, and the resulting SFT is not necessarily strongly aperiodic (in fact, if \(G\) is infinite this SFT's is necessarily not strongly aperiodic).  The ideas presented here are similar to ideas that date back to Gromov's paper~\cite{gromov}; see Coornaert and Papadopoulos~\cite{MR1222644,MR1878587} for a detailed treatment.

A key to our construction is {\it Cannon's Shortlex FSA}. Fixing an arbitrary order on  \(S\) induces a lexicographic order on the finite words \((S \cup S^{-1})^{*}\).  A path in the Cayley graph is called \em shortlex \em if and only if it is a geodesic, and is first in the lexicographic order among all geodesics with the given endpoints. It is clear that for each element \(g \in G\) there is a unique shortlex geodesic from the identity \(e\) to \(g\) (and this is true for any finitely generated group). It is far less obvious how calculate shortlex representatives. In fact, shortlex representatives in a group \(G\) can be calculated if and only if  \(G\) has a decidable word problem, see~\cite{book}.  

A hyperbolic group \(G\) admits an FSA,\footnote{For those unfamiliar, an FSA (finite state automaton) is a Turing machine without a tape; the ``memory'' is contained is finitely many states.}
described in~\cite{book},
that accepts a word if and only if it is shortlex representative. Denoting the states for the FSA as 
\[
\left\{
s_{1},\dots,s_{n}
\right\}
\]  
We can now label each group element with the following three labels:
\[
\left(
P(g), \ \mathrm{dist}(e,g), \ s_{i(g)}
\right)
\]
Defined as follows:
\begin{enumerate}
\item \(P(g)\) is the generator that points towards \(e\) on the shortlex geodesic to \(g\), in other words, the product \(gP(g)\) is the last group element on the shortlex geodesic before arriving at \(g\). \(P(e)\) is not defined. \(P\) is called the \em parent function\em.

\item \(\mathrm{dist}(e,g)\) is the distance from \(e\) to \(g\).  

\item \(s_{i(g)}\) is the state of the FSA at \(g\) (so \(i(g) \in \{1,\dots,n\}\)).
\end{enumerate}
Of course, \(\mathrm{dist}(e,g)\) cannot be used as a label for an SFT since it takes on infinitely many values. We replace it with the function that describes the \em difference \em of the distance to the origin between \(g\) and its neighbors, that is, for each \(a \in S\), we define
\[
\textrm{\dh}(g)(a) := 
\mathrm{dist}(e,g) -  \mathrm{dist}(e,ga)
\]
Since \(g\) and \(ga\) are neighbors this function can only take the values \(\pm 1\) or \(0\). Thus \(\textrm{\dh}\) is a function
\[
\textrm{\dh}: G \to
\{-1,0,1\}^{S}
\]
This completes the labeling on \(g\) for each group element. This labeling is \em not \em a configuration in any SFT but rather a blueprint for constructing \(\Omega_{S}\), as we now explain.

\bigskip\noindent
We now define a full shift of \(G\) with labels
\begin{equation}
\label{labels_shellings}
(\textrm{\dh}(g) ,s_{i(g)},P(g))
\in
\{-1,0,1\}^{S} \times \{s_1,\dots,s_{n}\} \times S
\end{equation}
Let \(\Omega_{S}\) be the subshift consisting of all configuration satisfying the following condition:
\begin{center}
The labeling of any ball, of any radius, \\ 
coincides with labeling above on some ball not containing \(e\).
\end{center}
It is not clear that:
\begin{enumerate}
\item This is an SFT. The requirement above says ``balls of any radius'', and we must ensure that this can be enforced by considering balls of a fixed radius. This is one of the many places that hyperbolicity is used.
\item This SFT is not empty.
\end{enumerate}
The second point is actually not too hard. Since we can label balls of arbitrary radius, a diagonalization argument shows that there is a labeling of all of \(G\), that is, the SFT is not empty.

\begin{remark}
\label{remark_notSA}
It is known that \(\Omega_{S}\) is necessarily \em not \em SA. More needs to be done.
\end{remark}

\begin{remark}
We give an example of a parent function \(P\) defined on \(F_{2}\), the free group on  generators \(a\) and \(b\). The Cayley graph of \(F_{2}\), a regular 4-tree, is shown in Figure~\ref{figure:papas}. 
At each vertex there are two gray edges (corresponding to \(a\) and \(a^{-1}\)) and  two black edges (corresponding to \(b\) and \(b^{-1}\)). Each edge is marked with an arrow that points from \(g\) to \(P(g)\).
As the figure suggests, \(P(g)\) can be \em any \em neighbor of \(g\); the only rule is that at any \(g\) we see exactly one triangle pointing ``out'' (towards \(P(g)\)) and three pointing in (from \(P^{-1}(g)\)). Picking any \(g \in F_{2}\) and following the arrows we arrive at a point on \(\partial_{\infty} F_{2}\), and this point is independent of choice of \(g\). 
\begin{figure}
\centerline{{\includegraphics[width=4.5in]{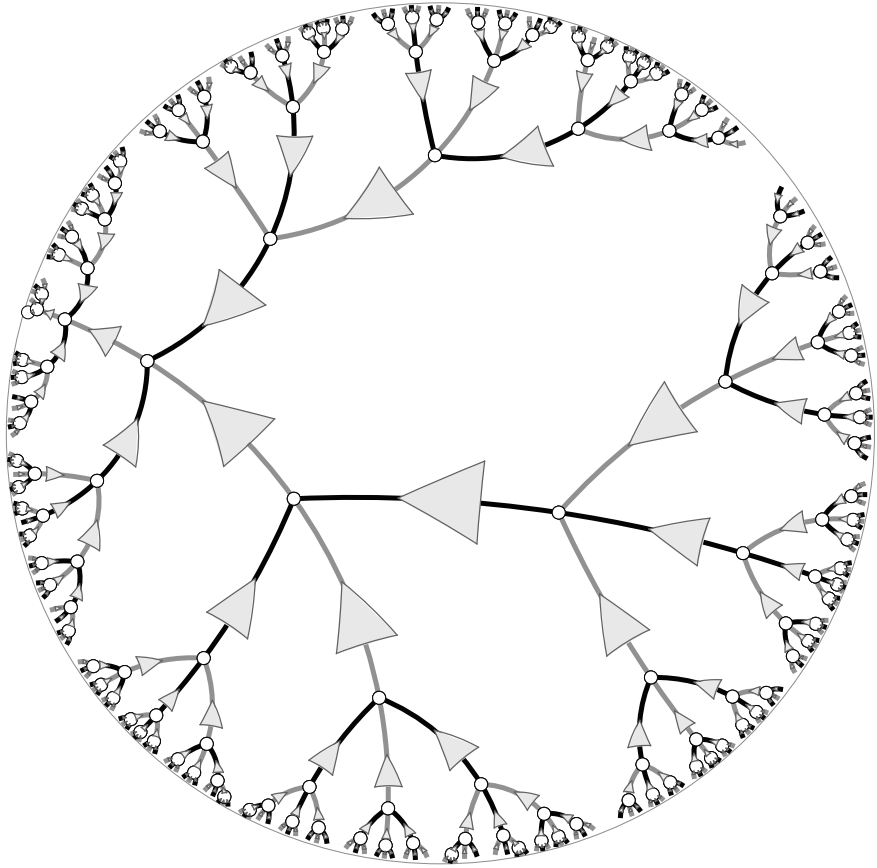}}}
\caption{The function \(P\)}
\label{figure:papas}
\end{figure}
\end{remark}

\section{Relation to the Boundary and Weak Aperiodicity}
\label{section:boundary_WA}

We are assuming that the reader is familiar with \(\partial_{\infty} G\), the \em boundary at infinity \em of \(G\).  

Given any configuration \(\omega \in \Omega_{S}\) and any \(g \in G\), we construct the path
\[
g, Pg, P^{2}g,\dots
\]
It follows from the definition of \(P\) that this path is a geodesic ray and hence defines a point at infinity. Hyperbolicity  implies that this point is independent of \(g\), and we denote it by \(\xi\) (of course, \(\xi(\omega)\) depends on \(\omega\)). The association \(\omega \mapsto \xi\) defines a function
\[
\Omega_{S} \to \partial_{\infty} G
\]
It is clear that this function is compatible with the actions of \(G\) on \(\Omega_{S}\) and on \(\partial_{\infty} G\). This means that the action of \(G\) on \(\partial_{\infty} G\) is a \em factor \em of an SFT, which motivated Gromov to study them.

We would like to exploit the function \(\Omega_{S} \to \partial_{\infty} G\) differently: we will use it to show that \(\Omega_{S}\) is weakly aperiodic. The map \(\Omega_{S} \to \partial_{\infty} G\) shows that for any \(\omega \in \Omega_{S}\) we have
\[
\stab[\omega] < \stab[\xi(\omega)]
\]
It is known that \(\stab[\xi]\) is virtually cyclic for any \(\xi \in \partial_{\infty} G\).\footnote{To see that \(\stab[\xi]\) is virtually cyclic use the following: (1) every non-torsion element of \(G\) fixes exactly two points on the boundary; (2) if two  elements share one fix points they share both fixed points; and (3) the elements that fix the same two points form a virtually cyclic group (this allows for a finite group, as finite groups are virtually trivial). See Sections~8.1 and~8.2 Gromov~\cite{gromov} and the proof of Proposition~III.\(\Gamma\).3.20 (Page~467) of~\cite{BridsonHaefliger}.}
We conclude that \(\stab[\omega]\) is virtually cyclic (for any \(\omega \in \Omega_{S}\)). By assumption \(G\) is one-ended, and hence \(G\) is not virtually cyclic. This shows that \(\stab[\omega]\) has infinite index, in other words, \(\Omega_{S}\) is weakly aperiodic.

\section{Horospheres}
\label{section:horospheres}

Fix a configuration \(\omega \in \Omega_{S}\).  We may ``integrate'' \dh\ (notation as in~~\eqref{labels_shellings} above) to get a function \(h:G \to \mathbb{Z}\) (defined up-to an additive constant).  The level sets of \(h\) are called {\it horospheres}. 

Any \(g \in \stab[\omega]\) preserves \dh\ and hence \(g\) preserves \(h\) up-to an additive constant, that is, there exist an integer \(C\) so that for any \(x \in G\) we have
\begin{equation}
h(g \cdot x) = h(x) + C
\end{equation}
If \(g\) is a torsion element then \(C=0\). Conversly, if \(g\) has infinite order then \(C \neq 0\) (Lemma~9.1  of~\cite{cohen_goodman-strauss_rieck_2021}). This simple observation will prove quite useful.

Since a hyperbolic group has finitely many conjugacy classes of torsion, getting rid of the torsion in the stabilizers is very easy (Proposition~3.3 of~\cite{cohen_goodman-strauss_rieck_2021}). Thus we may assume that \(g\) has infinite order and \(C \neq 0\).

\section{The Divergence Graph}

Fix \(\omega \in \Omega_{S}\).
As we discussed in Section~\ref{section:boundary_WA}, for every \(g \in G\) the sequence 
\begin{equation}
\label{past_geodesic}
g, Pg, P^{2}g,\dots
\end{equation}
converges to \(\xi(\omega) \in \partial_{\infty} G\). 
For convenience we describe \(P\) as moving ``down'' towards \(\xi\).

We now move up from  \(g \in G\):
\[
g, P^{-1}g, P^{-2}g,\dots
\]
This defines sets that move away from \(g\); unlike the downwards path that limits on \(\xi\), the limit of the sets \(P^{-n}g\) is more complex. It \em does \em depend on \(g\), and will often be an uncountable subset of \(\partial_{\infty} G\).  The union of this sets is called the {\it future cone} of \(g\), namely,
\[
P^{-*}g = \bigcup_{n=0}^{\infty} P^{-n}g
\]
\begin{remark}[growth rate]
\label{remark:growth_rate}
The hyperbolic group \(G\) has a well-defined growth rate which we will denote by \(\lambda > 0\).  Given a configuration \(\omega \in \Omega_{S}\), the future cone of each \(g \in G\) is completely determined by the state of the FSA at \(g\), denoted \(s_{i(g)}\) in~\eqref{labels_shellings} above (this was analyzed in detail in~\cite{MR4015648}). In particular, \(s_{i(g)}\) determines the growth rate of the future cone, and we call it the \em growth rate of the state\em. In what follows we only consider states whose future cone has growth rate \(\lambda\) (for a precise discussion see Definition~6.3 of~\cite{cohen_goodman-strauss_rieck_2021}).  We denote that set of all elements of \(G\) satisfying this condition \(G^{+}\). 
Note that each \(g \in G^{+}\) accumulates to an uncountable subset of \(\partial_{\infty} G\). 
By Proposition~6.5 of~\cite{cohen_goodman-strauss_rieck_2021} any \(2\delta\)-ball, in any configuration \(\omega \in \Omega_{S}\), contains a point of \(G^{+}\). 
\end{remark}

Back to our discussion, fix a configuration \(\omega \in \Omega_{S}\).  

Let  \(h\) be as in Section~\ref{section:horospheres}.
For \(i \in \mathbb{Z}\) set \(H_{i} := h^{-1}(i)\) and \(H^{+}_{i} := G^{+} \cap H_{i}\).  On \(H^{+}_{i}\) we define the {\it divergence graph} as follows:
\begin{enumerate}
\item[\bf v] The vertices of the divergence graph are the vertices of \(H^{+}_{i}\).
\item[\bf e] Two vertices \(g_{1}, g_{2} \in H^{+}_{i}\) are connected by an edge if and only if their futures remain a bounded distance apart. In other words, for some \(C>0\), and for each integer \(n \geq 1\), there are \(v_{1} \in P^{-n}g_{1}\) and \(v_{2} \in P^{-n}g_{2}\) with \(\mathrm{dist}(v_{1},v_{2}) < C\). 
\end{enumerate}

The following holds (see Lemma~7.4  of~\cite{cohen_goodman-strauss_rieck_2021} and its proof):
\begin{itemize}
\item The union of the limit sets of the future cones of all of \(H^{+}_{i}\) is \(\partial_{\infty} G - \{\xi\}\).

\item \(g_{1}, g_{2} \in H^{+}_{i}\) are connected by a divergence graph edge if and only if the limit sets of their future cones intersect.

\item We may therefore describe the divergence graph as a discrete approximation of \(\partial_{\infty} G\).

\item Most importantly, \em the divergence graph is connected\em. This reflects the fact that \(\partial_{\infty} G - \xi\) is connected (the Cut Point Conjecture, proved by Swarup~\cite{MR1412948}). 
\end{itemize}

The plan is now as follows. Any infinite order element \(\phi \in \stab[\omega]\) will translate the levels of \(h\) by \(C \neq 0\), as explained in Section~\ref{section:horospheres}. It is our goal to enhance \(\Omega_{S}\) by associating an integer \(\Delta(g)\) to each \(g \in G^{+}\) in a way that cannot be periodic; the new SFT will no longer have \(\phi\) in its stabilizer.

We enhance the labels of \(\Omega_{S}\) (compare this with the labels presented in~\eqref{labels_shellings})

\begin{equation}
\label{labels_populated_shellings}
(\textrm{\dh}(g) ,s_{i(g)},P(g),\wp(g),\Delta(g),m(g))
\end{equation}

The enhanced SFT is called \em populated shelling\em, denoted \(\Omega_{P}\). The name comes from the fact that \(\wp\) defines a ``population'' of ``villagers'' on each ``village'' \(g \in G\) (or, if the reader prefers, any village \(v \in G^{+}\), since the population of any \(g \in G - G^{+}\) is zero anyway).

\section{\(\wp\) and \(m\)}
\label{section:PandM}

The number of ``villagers'' defines the function \(\wp\), with \(\wp(g)\) being the population at \(g\):
\[
\wp: G^{+} \to \mathbb{Z}_{\geq 1}
\]
This function is required to be bounded, with the population bound \(N\) fixed in advanced.

Clearly, we need something that will help us relate \(\wp\) to the geometry of the group, for otherwise the population values will be arbitrary numbers. This is \(m\), which stands for \em matching\em. More precisely, it is \em parent-child matching\em. 

Each ``villager'' in \(v \in H^{+}_{i}\) has \(q^{\Delta_{i}}\) children (\(q^{\Delta_{i}}\) will be described in the next section, for now just take it to be ``some number''). We list the children as follows: 
\[
(v,j,k)
\]
Where here \(v\) as the village, \(1 \leq j \leq \wp(v)\) is the villager, and \(1 \leq k \leq q^{\Delta_{i}}\) is the child. The function \(m\) ``places'' this child as a villager in \(H^{+}_{i+1}\):
\[
m(v,j,k) = (u,l)
\]
Here, \(u \in H^{+}_{I+1}\) is a village and \(1 \leq l \leq \wp(u)\) is a villager.

The geometry of the group comes to play when we force the child to be placed not too far from the parent.  The precise condition is the following: the child may take up to 3 steps on the divergence graph of \(H^{+}_{i}\), and then move one step up. A succinct description is this: the child of a villager in \(v \in H^{+}_{i}\) is a villager in \(u  \in H^{+}_{i+1}\) with 
\[
\mathrm{DivDist}(v,Pu) \leq 3
\]
Naturally, \(\mathrm{DivDist}\) denotes that divergence graph distance. 

\begin{remark}
The reader probably finds the constant 3 rather arbitrary (not to say mysterious). It comes from an application of a theorem in graph theory, which states that if a graph is connected, then its cube 
admits a Hamiltonian path between nay two vertices.\footnote{The \em cube \em of a graph is obtained by adding an edge between any two vertices of distance at most 3.} 
This is then used to construct a ``translation-like \(\mathbb{Z}\) action'' (in the sense of Seward~\cite{MR3158775}) whose {\it defect}  is 3, and this is the origin of the constant. 
\end{remark}

\section{Straying Away and Coming Back Home}

It would be quite natural to worry that we are too loose with the geometry of the group here. We populate the group with the goal of considering the population after arbitrarily many generations, and descendants may stray 3 divergence-graph steps each generation. 


Indeed, after many generation, a descendent of the villager \((v,j)\) may be in a village which is very far from the future cone \(P^{-*}v\). It may be worth emphasizing that we are using two distinct ``futures'' here, the future cone \(P^{-*}v\) which is the collection of \em villages \em \(u\) for which
\[
v = P^{n}u
\]
for some \(n\). On the other hand, there are \em villagers \em that are descendants of  \((v,j)\), and may stray away from \(P^{-*}v\).  

Hyperbolic geometry to the rescue. As described above, let \((u,l)\) be a descendant of \((v,j)\), after, say, \(n\) generations. We use the notation \(v=v_{0},v_1, \dots, v_{n}=u\) for the villages so that \((u,l)\) is a descendant of \((v_{n-1},l_{n-1})\) (for some \(l_{n-1}\)), \((v_{n-1},l_{n-1})\) is a descendant of \((v_{n-2},l_{n-2})\) (for some \(l_{n-2}\)), and so one. It is not hard to show that an edge of a divergence graph connects vertices of Cayley distance at most \(2\delta\), and so the apple doesn't fall too far from the tree:
\[
\mathrm{CayDist}(v_{i},v_{i+1}) \leq 6\delta + 1
\]
On the other hand, because \(v_{i} \in H^{+}_{i}\) and \(v_{j} \in H^{+}_{j}\), we get that
\[
\mathrm{CayDist}(v_{i},v_{j}) \geq |i-j|
\]
This produces a \em quasi-geodesic \em that can be compared with the geodesic
\[
v_{n},Pv_{n},\dots,P^{n}v_{n}
\]
It is a feature of hyperbolic geometry that quasi-geodesics remain a bounded distance away from geodesics, which means that, for some fixed \(R>0\) we have:
\[
\mathrm{CayDist}(v_{0},P^{n}v_{n}) \leq R
\]
A precise statement is given in Lemma~9.2 of~\cite{cohen_goodman-strauss_rieck_2021}. This allows for sufficient control over the population growth, since it shows that all the descendants of villagers in \(S \subset H^{+}_{i}\), where \(S\) is any finite set, live in the future cone of \(\mathcal{N}(S)\), where here \(\mathcal{N}(S)\) is the 3-neighborhood of \(S\) in the divergence graph on \(H^{+}_{i}\).

\section{\(\Delta\)}
\label{section:Delta}

We finally describe \(\Delta\), focusing on \(G^{+}\); \(\Delta\) is extended to \(G\) by setting it to be zero on \(G - G^{+}\).
 It is \(\Delta\) that will ultimately be responsible for aperiodicity.

For each \(g \in G^{+}\), the growth rate of the population at \(g\) is controlled by \(\Delta(g)\). To be precise, every villager at \(g\) has exactly \(q^{\Delta(g)}\) children (for \(q\) to be decided momentarily). 
We use, intentionally, a number \(q\) which is \em not \em compatible with \(\lambda\):  \(q\) is an integer, which we may take to be either 2 or 3, so that
\begin{equation}
\frac{\log(q)}{\log(\lambda)}
\not\in \mathbb{Q}
\end{equation}
The function
\[
\Delta: G^{+} \to 
\left\{
\left\lfloor \log_{q}(\lambda) \right\rfloor, \left \lceil \log_{q}(\lambda)\right \rceil
\right\}
\]
is required to satisfy the following conditions:
\begin{enumerate}

\item (And this is key) \(\Delta\) is constant along levels of \(h\), that is, for each integer \(i\), \(\Delta|_{H^{+}_{i}}\) is constant. We denote this value \(\Delta_{i}\).  This defines the sequence 
\[
\left(\Delta_{i}\right)_{i \in \mathbb{Z}}
\]

\item The sequence \(\Delta_{i}\) approximates \(\lambda\) (in a sense made precise in Corollary~9.3 of~\cite{cohen_goodman-strauss_rieck_2021}).  

\item The condition above, and the incompatibility of \(\lambda\) and \(q\), guarantee that \(\Delta_{i}\) is not periodic  (Corollary~9.4 of~\cite{cohen_goodman-strauss_rieck_2021}).  
\end{enumerate}

We saw in Section~\ref{section:horospheres} that any infinite order elements in the stabilizer of a configuration \(\omega \in \Omega_{S}\) must translate the levels \(\{h=i\}\) by a non zero amount. Once we enhance \(\Omega_{S}\) by populating \(G^{+}\), the function \(\Delta\) will not be invariant under such an element. This shows that the there is no infinite-order element in the stabilizer of any configuration in \(\Omega_{P}\), as desired. 

\begin{remark}[the role of the (necessary!) assumption of one-endedness of \(G\)] 
The question of existence of SA SFT is irrelevant for zero-ended groups (that is, finite groups), where the answer is always ``yes'', as well as two-ended groups (that is, virtually-\(\mathbb{Z}\)'s), where the answer is always ``no''. So we ignore these groups in this remark and consider only one-ended and infinitely-ended groups. Our main result is that an SA SFT exists \em only \em for the former. This begs the question: where exactly was the ``one-endedness'' assumption used? The only place is {\it imposing that \(\Delta\) be constant along each \(H^{+}_{i}\)}. This must be enforced via local rules, as part of the SFT. What allows us to do this is \em connectivity of the divergence graph\em. This is the one and only place where the assumption is used, as connectivity of the divergence graph is equivalent to \(G\) being one-end. Swarup's resolution of the cut-point conjecture~\cite{MR1412948} plays a key role here; it states that \(\partial_{\infty} G \setminus \{\xi\}\) is connect (for any \(\xi \in \partial_{\infty} G\)).
\end{remark}

\section{One Last Issue}

It is not the goal of this paper to give a complete proof of theorem~\ref{MainThm}. Our goal is to explain some of the elements that go into the proof in a way that would facilitate its reading. However, it is hard to ignore the fact that we have not addressed the following question:
\begin{center}
Does a populated shelling even exist?
\end{center}
This should not be taken lightly as it is quite possible that we defined the empty subshift here. In fact, if the population bound is too small this is probably the case. Proposition~8.5 of~\cite{cohen_goodman-strauss_rieck_2021} shows that \(\Omega_{P}\) described above is indeed an SFT, and Proposition~9.5 shows that no configuration has an infinite order element in its stabilizer, but neither addresses existence of a configuration.

Much of the work in~\cite{cohen_goodman-strauss_rieck_2021} is devoted to Proposition~8.12, showing that (for an appropriately chosen population bound) \(\Omega_{P}\) is indeed not empty. This is the most technical and longest part of the proof, and here is my attempt at explaining the idea. We proceed in three steps:

\begin{enumerate}
\item[\bf Level:] The first step is populating each level \(H^{+}_{i}\). To discuss that we need to dig a little deeper into the maximal growth states of the FSA.\footnote{ 
There are two distinct notions of growth rate at play, growth rate of the population, controlled by \(\Delta_{i}\), and growth rate of the group elements, which is \(\lambda\); here we discuss the latter. Having two incompatible growth rates is a little confusing but it is the very thing the leads to strong aperiodicity.}
In Remark~\ref{remark:growth_rate} we explained that the future of some states must have growth rate \(\lambda\) (the growth rate of the group itself) and denoted the set of vertices that have this growth rate by \(G^{+}\). From that point we concentrated on \(G^{+}\) and on \(H^{+}_{i} := G^{+} \cap H_{i}\). In fact more is true; although the transitions of the FSA need not satisfy the assumptions of the Perron--Frobenius Theorem, it is possible to associate to them a measure that behaves just like a Perron--Frobenius eigenvector; see~\cite{MR4015648} or Section~6 of~\cite{cohen_goodman-strauss_rieck_2021}. Now in the first step we populate each level \(H^{+}_{i}\) so that the population of \em any \em finite subset approximates its total measure (it is not possible to get this to be exact; the ratio of the measures of distinct states is usually irrational). This is Lemma~8.7 of~\cite{cohen_goodman-strauss_rieck_2021} (for one level). The average ratio of population to total measure (called the population density) can be chosen freely, within a reasonable range (Definition~8.6 and Lemma~8.7).

\item[\bf Levels:] We apply this to all of \(G\) (one level at a time, and still without matching levels). This is given in Corollary~8.9 of~\cite{cohen_goodman-strauss_rieck_2021}. In order to be able to match parents and children, we ensure that the sequence of densities behaves well (grows when small, shrinks when big). The precise description is given in Definition~8.10.

\item[\bf Matching:] Having populated the group as described above, we apply the Hall Matching Theorem to prove existence of a matching function \(m\) as required.  This is Proposition~8.11 of~\cite{cohen_goodman-strauss_rieck_2021}.

\item[\bf Ta da!] That's all, folks.

\end{enumerate}

\end{document}